# The Hamiltonian properties of rectangular meshes with at most two faulty nodes

Yingtai Xie

Chengdu university

**Abstract.**

In this paper we consider the Hamilton cycle problem in the rectangular meshes with at most two faulty nodes. We prove that this problem is solvable in polynomial time with a corresponding algorithm. We provided an entirely new approach to this problem being different from the method early used on this problem.

1. **introduction.**

In the design of parallel algorithms, In order to efficiently implement a linear-array-based or ring-based parallel algorithm on a mesh-connected parallel computer, one needs to find a long path/cycle within the mesh network. Chen et al. [1] presented a fast sequential algorithm as well as its parallelized version for constructing Hamiltonian paths between two nodes in a mesh. Collins and Krompart [2] and Stoyan and Strehl [3] studied the number of Hamiltonian paths and Hamiltonian cycles in a mesh, respectively.

There may be failing nodes in a huge number of mesh-connected computer systems. Once the failing nodes are identified, the next work is to identify a long fault-free path/cycle. It is a Hamilton cycle/path problem in a 2D-mesh with faulty nodes. Unfortunately, it is NP-complete for arbitrary faulty nodes [Itai et al. 4] so ther is no polynomial-time algorithm for arbitrary faulty nodes unless P=NPcom.

Hedetniemi et al. [5] proposed a ring embedding method in a mesh with a single faulty node. Kim and Yoon [6] found a necessary and sufficient condition for a Hamiltonian $m \times n$ mesh ($m$ and $n$ are multiples of four) in the presence of two faulty nodes; they also devised an algorithm for constructing such a Hamiltonian cycle, if there is, by partitioning the original mesh into $4 \times 4$ submeshes, building an appropriate fault-free Hamiltonian path for each submesh, and merging these short paths to yield a fault-free Hamiltonian cycle of the original mesh. Xiaofan Yang et al [7] describe a systematic method for constructing such a Hamiltonian cycle, if there is; their method is different from the one used by Kim and Yoon in that the original mesh is partitioned into at most five submeshes instead of a large number of $4 \times 4$ submeshes. For some related work the reader is referred to [8,9,10,11].

In this paper we provided a more concise necessary and sufficient condition and more explicit algorithm for having a Hamilton cycle in odd size $m \times n$ mesh with one faulty node and in even $m \times n$ mesh with two faulty nodes.

2. **Preliminaries**

A $m \times n$ (rectangular) mesh is isomorphic to a square grid graph $R(m,n)$ whose vertex set is

$$V(R(m,n)) = \{v : 0 \le v_x \le m-1; 0 \le v_y \le n-1\}$$

and an edge $(u,v) \in E(R(m,n))$ if and only if $|u_x - v_x| + |u_y - v_y| = 1$

A $m \times n$ (rectangular) mesh with faulty nodes is isomorphic to a square grid graph $R(m,n)$-{u} with holes formed after removing some vertices {u}.

We say that vertex $v$ is even(white dots throughout this paper) if $v_x + v_y \equiv 0$ (mod 2); otherwise, $v$ is odd(black dots). It is immediate that $R(m,n)$ is bipartite, with the edges connecting even and odd vertices. We say that $R(m,n)$ is even-sized if $m \times n$ is even and it is odd-sized if $m \times n$ is odd.

2. **Hamilton cycles and paths in rectangular .**

We first examine the existence of Hamilton cycles and paths in $R(m,n)$ for a given array $m,n$. To avoid trivialities, it is necessarily that $m,n>1$. We will prove that if $R(m,n)$ is even-sized then there is a Hamilton cycle. If $R(m,n)$ is odd-sized then there is no Hamilton cycle in it, because it is bipartite, but it has a Hamilton path.

An unit square is a simple cycle with minimal length. Let $C_1$ and $C_2$ be two disjoint cycles. An unit square $C_f$ is called a separant of $C_1$ and $C_2$ if the two paralleled edges of $C_f$ share with $C_1$ and $C_2$ respectively. (Fig 2.1 (a))

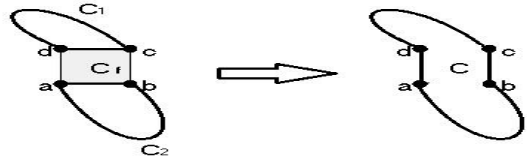

(a) A separant $C_f$ of $C_1$ and $C_2$    (b) $C=(C_1\cup C_2)\oplus C_f$

Fig.2.1    Cycle-merging operation

***Cycle-merging operation***：
$$C=(C_1\cup C_2)\oplus C_f \qquad (2.1)$$

$C$ is a cycle containing all of the vertices of $C_1$ and $C_2$, formed by removing edges ab and cd, which are shared with $C_f$, from $C_1$ and $C_2$ and adding two edges bc and cd of $C_f$.

Let a path $P_1$ and a cycle $C_1$ be disjoint, an unit square $C_f$ is called a ***separant*** of $P_1$ and $C_1$ if the two parallel edges of $C_f$ share with $P_1$ and $C_1$ respectively. (Fig 2.2 (a))

***Cycle-path-merging operation***：
$$P=(P_1\cup C_1)\oplus C_f \qquad (2.2)$$

(refer to Cycle-merging operation) $P$ is a path containing all of the vertices of $P_1$ and $C_1$.

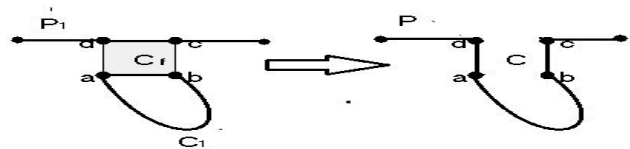

(a)        (b)

Fig.2.2    Cycle-path-merging operation

*A connected subgraph* $C$ *of* $G$ *is a cycle* if $d_C(v)=2$ for each $v$ on $C$ and it is a Hamilton cycle if $V(C)=V(G)$.

*A connected subgraph* $P(v_1,v_2,...,v_k)$ *of* $G$ *is a path* if for two special vertices, called end-points of $P$, $v_1, v_k$, $d_p(v_1)=d_p(v_k)=1$; All the other points, called inner points, $d_p(v_i)=2, 1<i<k$.

A ***2-factor*** $F(C_1,C_2,....C_m)$ of a graph $G$ is a set of $m$ disjoint cycles that span $G$. when $m=1$, it is a Hamilton cycle. $F$ is a 2-factor of $G$ if and only if $d_F(v)=2$ for each $v\in V(G)$.

A **[1,2]-*factor*** $D(P_0,C_1,...C_m)$ of a graph $G$ is set of a path and $m$ disjoint cycles that span $G$. when $m=0$, it is a Hamilton path. $D$ is a [1,2]-factor of $G$ if and only if for two special points $s,t$ (the end of $P_0$), $d_D(s)=d_D(t)=1$, all other points $d_D(v)=2$ for each $v\in V(G)$.

We say that two disjoint subgraph $G_1$ and $G_2$ are **neighbors** if there is a edge $(u,v)\in E(G)$ such that $u\in V(G_1)$ and $v\in V(G_2)$.

**Observation 2.1.** let $F(C_1,C_2,....C_m)$ be 2-factor of a connected graph $G$ then each cycle

$C_i$ have at least a neighbor $C_j$ $(i \neq j)$ .

**Observation 2.2.** let $D(P_0,C_1,...C_m)$ be a [1,2] factor of a connected graph $G$ then path $P_0$ have at least a neighbor $C_k$ and each cycle $C_i$ have at least a neighbor $C_j$ $(i \neq j)$.

It is immediate by inductively doing cycle-merging operation and cycle-path-merging operation that following:

**Lemma 2.1** Let $F(C_1,C_2,....C_m)$ is a 2-factor of a grid graph $G$. If there is a separant of any two neighbors $\{C_i,C_j\}$ then a Hamilton cycle can be constructed in $G$.

**Lemma 2.2** Let $D(P_0,C_1,C_2,....C_m)$ is a [1,2]-factor of a grid graph $G$. If there is a separant of any two neighbors $\{C_i,C_j\}$ or $\{P_0,C_k\}$ then a Hamilton path can be constructed in $G$.

We have the following theorem.

**Theorem 2.1** $R(m,n)$ has a Hamilton cycle if and only if it is even-sized and it has a Hamilton path if it is odd-sized.

**Proof** : Let $R(m,n)$ be even-sized, without loss of generality, let $m$ be even. A strip is a rectangular graph with minimum dimension 2. The strip with corners $a,b,c,d$ ( $a$ is adjacent to $b$ and $d$ ; $b$ is adjacent to $c$ ; $c$ is adjacent to $d$ ) is denoted $S(a,b;c,d)$.(Fig.2.3(a)). The perimeter of a strip is a cycle. The $\frac{m}{2}$ strips with corners $a_i=(2i,0)$, $b_i=(2i+1,0)$ and $c_i=(2i+1,n-1)$, $d_i=(2i,n-1)$ for $i=0,1,2,...\frac{m}{2}$, formed a 2-factor of $R(m,n)$ and an unit square between $S_i$ and $S_{i+1}$ is a separant of them. A Hamilton cycle can be constructed by lemma 2.1.(refer to Fig.2.3 (b),(c))

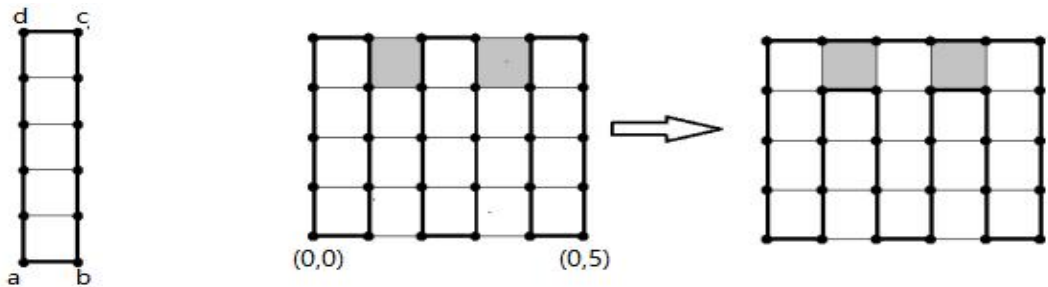

(a)   a strip          (b)   2-factor of $R(6,5)$          (c)   A Hamilton cycle of $R(6,5)$

Fig.2.3 Construct a Hamilton cycle in $R(6,5)$

Let $R(m,n)$ be odd-sized then there are a [1,2]-factor formed by $\frac{1}{2}(m-1)$ strips and a path and respective separants, a Hanilton path can be constructed by lemma 2.2 (refer to fig.2.4)

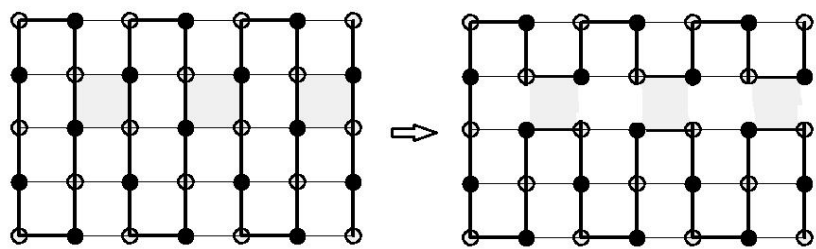

(a) A [1,2]-factor of $R(7,5)$          (b) A Hamilton path of $R(7,5)$

Fig.2.4 Construct a Hamilton path in $R(7,5)$

## 3. The Hamilton Cycle in rectangular grid graph with missing points.

A mesh with faulty nodes isomorphic to a rectangular grid graph with missing points.

If there is a 2-factor in a rectangular grid graph with missing points then there is a Hamilton cycle in it by lemma 2.1.

What we should to do is that to find the 2-factor in the rectangular grid graph with at most two missing points. And then, to construct a Hamilton cycle in it.

We call $M$ **2-limited spanning subgraph** of $G$ if $V(M) = V(G)$ and $1 \le d_M(v) \le 2$ for each $v \in V(G)$. A $v \in V(G)$ is called **filled** if $d_M(v) = 2$, otherwise **unfilled**. Let $\sigma(M)$ be the number of unfilled vertices in $M$ then $M$ is 2-factor of $G$ if $\sigma(M) = 0$.

we ascribe edges of $M$ red and edges of $E(G) - M$ blue.

**Definition 3.1** A path $P$ is called **M-alternating path** if its edges are alternately blue and red.

**Definition 3.11** A $M$-alternating path
$$P = (v_0, v_1, \ldots v_{k-1}, v_k) \qquad (3.11)$$
meeting the condition that its end points $v_0$ and $v_k$ both are unfilled and both of the first edge $v_0 v_1$ and the terminal edge $v_{k-1} v_k$ are blue, is called **M- augmenting path.**

If $P$ is a $M$-augmenting path then $|P|$ must be odd and $v_0$ and $v_k$ have different colour.

**Operation** $\oplus$ : Let $M_1$ and $M_2$ be two sets
$$M_1 \oplus M_2 = (M_1 \cup M_2) - (M_1 \cap M_2)$$

**Lemma 3.1** let $\sigma(M) \ne 0$ and $P$ is a M-augmenting path with end points $v_0$ and $v_k$ then
$$M_1 = M \oplus P \qquad (3.1)$$
Is also a 2-limited subgraph of $G$ and $\sigma(M_1) = \sigma(M) - 2$.

**proof**: We ascribe the edges of $M$ red and the edges of $E(G) - M$ blue. $M_1$ is formed by deleting all the red edges of $P$ from $M$ and adjoin to it all blue edges of $P$. It is clear only unfilled vertices $v_0$ and $v_k$ have each added one edge to become filled. Therefore $\sigma(M_1) = \sigma(M) - 2$. □

In simple terms, let $s$ and $t$ both be unfilled then a M-augmenting path $P$ from $s$ to $t$ convert them filled, and the M-augmenting pah $P$ itself be convert to a **cross path** of $M_1$

The following theorem is immediate.

**lemma 3.2** Let $M$ be a 2-limited spanning subgraph of grid $G$ and $\sigma(M) = 2r$ if there are $r$ disjoint $M$-augmenting paths then there is a 2-factor of $G$.

**Theorem 3.1** Let $R(m,n) - \{u\}$ is a odd rectangular grid graph with one missing point $u$, if $u$ is even then there is a Hamilton cycle in $R(m,n) - \{u\}$.

Proof: $R(m,n)$ is odd and $u$ is even. $R(m,n)$ have [1,2]-factor
$$D = C_1 \cup C_2 \cup \ldots \cup C_{(m-1)/2} \cup P_0$$

Let $x, y$ be two end-points of $P_0$ then $x, y$ are even, because $m \times n$ is odd. There

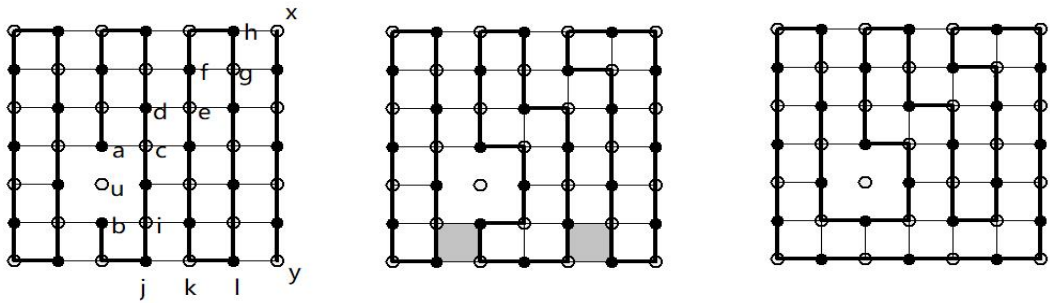

(a) $R(7,7) - \{u\}$    (b) A 2-factor in $R(7,7) - \{u\}$    (c) A Hamilton cycle in $R(7,7) - \{u\}$

Fig 3.1 A instance of theorem 3.1

are two vertices $a,b$ being of odd and being degree-1 in $M = D - \{u\}$. there are two M-augmenting path $P_1 = (a,c,d,e,f,g,h,x)$ and $P_2 = (b,i,j,k,l,y)$ (fig.3.1(a))
$$M_1 = (M \oplus P_1) \oplus P_2$$
Is a 2-factor including three cycles of $R(7,7) - \{u\}$. (Fig 3.1(b)). Do Cycle-merging operation for $M_1$ to construct a Hamilton cycle of $R(7,7) - \{u\}$. (Fig.3.1(c)). □

**Theorem 3.2** Let $R(m,n) - \{u,v\}$ is a even(let $m$ be even, $m,n \geq 4$) rectangular grid graph with one missing points $u,v$, meet the conditions:
(1) $u,v$ are different colors.
(2) If one of $u,v$ is a corner node of the original $R(m,n)$, then it is necessarily that the other is next to it..
Then there is a Hamilton cycle in $R(m,n) - \{u\}$.

Proof: $R(m,n)$ is even then $R(m,n)$ have 2-factor
$$C = C_1 \cup C_2 \cup ... \cup C_{m/2}$$

There are four vertices being of odd and being degree-1 in $M = D - \{u\}$. there are two M-augmenting path $P_1 = (a,c,d,e,f,g,h,x)$ and $P_2 = (b,i,j,k,l,y)$ (fig.3.1(a))
$$M_1 = (M \oplus P_1) \oplus P_2$$
Is a 2-factor including three cycles of $R(7,7) - \{u\}$. (Fig 3.1(b)). Do Cycle-merging

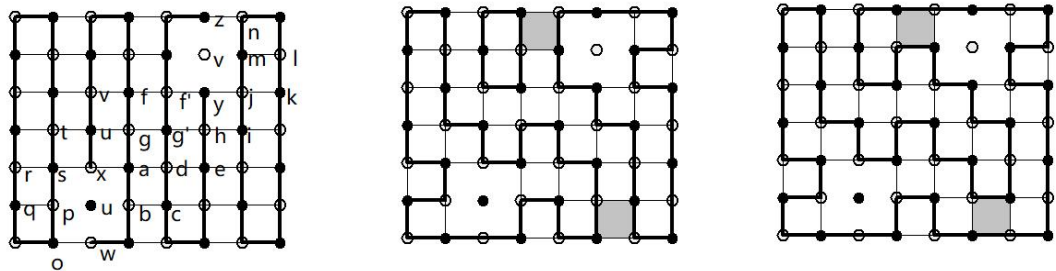

operation for $M_1$ to construct a Hamilton cycle of $R(7,7) - \{u\}$. (Fig.3.1(c)).

(a) $R(8,7) - \{u,v\}$     (b) A 2-factor in $R(8,7) - \{u,v\}$     (c) A Hamilton cycle in $R(8,7) - \{u,v\}$
Fig 3.2 A instance of theorem 3.2

# References

[1] S.D. Chen, H. Shen, R. Topor, An efficient algorithm for constructing Hamiltonian paths in meshes, Parallel Computing 28 (9) (2002)

[2] K.L. Collins, L.B. Krompart, The number of Hamiltonian paths in a rectangular grid, Discrete Mathematics 169 (1–3) (1997) 29–38.1293–1305.

[3] R. Stoyan, V. Strehl, Enumeration of Hamiltonian circuits in rectangular grids, Journal of Combinatorial Mathematics and Combinatorial Computing 21 (1996) 109–127.

[4] A. Itai, C. Papadimitriou, J.L. Szwarcfiter, Hamilton paths in grid graphs, SIAM Journal on Computing 11 (4) (1982) 676–686.

[5] S.M. Hedetniemi, S.T. Hedetniemi, P.J. Slater, Which grids are Hamiltonian? Congressus Numerantium 29 (1980) 511–524.

[6] J.S. Kim, S.H. Yoon, Embedding of rings in 2-D meshes and tori with faulty nodes, Journal of Systems Architecture 43 (9) (1997)643–654.

[7]Xiaofan Yang, Jun Luo, Shuangqing Li,Fault-tolerant Hamiltonicity in a class of faulty meshes, Applied Mathematics and Computation,Volume 182, Issue 2,2006,Pages 1696-1708.

[8] G.M. Megson, X. Liu, X. Yang, Fault-tolerant ring embedding in a honeycomb torus with nodes failures, Parallel Processing Letters
9 (4) (1999) 551–561.

[9] B. Parhami, An Introduction to Parallel Processing: Algorithms and Architectures, Plenum Press, 1999.

[10] K.W. Park, H.S. Lim, J.H. Park, H.C. Kim, Fault Hamiltonicity of meshes with two wraparound edges, in: Proceedings of the 10th
Annual International Conference on Computing and Combinatorics, Lecture Notes in Computer Science, 2004, pp. 412–421.

[11] H. Cho, L. Hsu, Ring embedding in faulty honeycomb rectangular torus, Information Processing Letters 84 (5) (2002) 277–284.